\documentclass[11pt]{article}

\if@twoside
\else
\fi
\usepackage[latin1]{inputenc}
\usepackage{makeidx}
\usepackage[french]{babel}
\usepackage{amsfonts}
\usepackage{amssymb}
\usepackage{amsthm}
\usepackage{newlfont}
\usepackage{showidx}
\usepackage[dvips]{graphicx}

\newtheorem{theo}{Th\'eor\`eme}
\newtheorem{prop}{Proposition}

\def\g{\mathfrak{g}}
\def\G{\Gamma}

\def\ov{\overline}
\def\tr{\mathrm{tr}}
\def\ad{\mathrm{ad}}

\newcommand{\dif}[2]{\frac{\mathrm{\displaystyle d}{\displaystyle #1}}{\mathrm{\displaystyle d}{\displaystyle #2}}}
\newcommand{\R}{\mathrm{ I\! R}}

\begin{document}

\bibliographystyle{alpha}
\author{Charles Torossian\footnote{UMR 8553 du CNRS, DMA-ENS, 45 rue d'Ulm 75230 Paris cedex 05}}
\title{ Sur la conjecture combinatoire de Kashiwara-Vergne}
\maketitle

{\noindent \bf{R\'esum\'e: }} Nous expliquons dans cet article un lien  entre les travaux
de Kontsevich sur la quantification formelle des
variétés de Poisson et la conjecture combinatoire de Kashiwara-Vergne.

{\noindent \bf{Abstract:} }We explain  a connection between the
combinatorial Kashiwara-Vergne conjecture
and the Kontsevich formula for quantization of Poisson manifolds.

\section{Introduction}
\subsection{La conjecture combinatoire de Kashiwara-Vergne}

Soit $\g$ une algèbre de Lie de dimension finie sur $\R$. Notons $\exp(x)$
l'image de $x\in \g$ par l'application exponentielle et nous noterons en général
$Z(x,y)$ la série de Campbell-Hausdorff définie par $Z(x,y)=\log(\exp(x)\exp(y))$.
 Dans tout ce qui suit nous travaillons
au niveau des séries formelles mais des arguments élémentaires montrent que toutes les
séries formelles que nous manipulons sont convergentes dans un voisinage proche de $0$.

La conjecture combinatoire de Kashiwara-Vergne \cite{KV} s'énonce de la manière suivante.
Il existe des séries $F(x,y)$ et $G(x,y)$ sur $\g\oplus\g$ sans terme constant
et à valeurs dans $\g$ telles que
l'on ait
\begin{equation}\label{1}
x+y-\log(\exp(x)\exp(y))=(\exp(\mathrm{ad}x)-1)F(x,y)+ (1-\exp(\mathrm{ad}y))G(x,y)
\end{equation}
et telle que l'identité suivante soit vérifiée
\begin{equation}\label{trace}{\mathrm tr}_{\g}(\mathrm{ad}x\partial_{x}F+\mathrm{ad}y\partial_{y}G)=
\frac{1}{2}{\mathrm tr}_{\g} (\frac{\mathrm{ad}x}{\exp(\mathrm{ad}x)-1}+\frac{\mathrm{ad}y}{\exp(\mathrm{ad}y)-1}
-\frac{\mathrm{ad}Z}{\exp(\mathrm{ad}Z)-1}-1).
\end{equation}

Cette conjecture combinatoire est démontrée dans l'article de Kashiwara-Vergne \cite{KV}
pour les algèbres de
 Lie résolubles. Elle fut démontrée pour $\mathrm{sl}(2)$ par Rouvière \cite{rou81}
et récemment par M. Vergne \cite{Ve}
 dans le cas quadratique (on rappelle que quadratique veut dire qu'il existe une forme bilinéaire
invariante et non dégénérée, mais on ne suppose pas que c'est la forme de Killing!). Il existe
une conjecture analogue dans le cas des espaces symétriques \cite{rou86}, \cite{rou90},
\cite{rou91}, \cite{rou94}.

Cette égalité sur les traces peut sembler étrange mais elle est une conséquence naturelle de
l'intégration par partie. Expliquons un peu tout ceci ce qui motivera le lecteur. Dans l'article
de Kashiwara-Vergne une idée de base est de considérer la déformation naturelle de  l'algèbre
de Lie $\g$ qui consiste à déformer le crochet $[x,y]$ en $t[x,y]$. Par exemple la série de
Campbell-Hausdorff est changée en $Z_{t}(x,y)=\frac{1}{t}Z(tx,ty)$.

On déduit facilement de l'équation (\ref{1}) et de la formule de la différentielle de l'application
exponentielle que l'on a une équation différentielle de la forme
\begin{equation}\label{eqdiff}
\frac{\partial}{\partial t}Z_{t}(x,y)=[x, F_{t}(x,y))]\cdot\partial_{x}Z_{t} +[y, G_{t}(x,y))]\cdot\partial_{y}Z_{t}
\end{equation}
où on a noté $F_{t}(x,y)=\frac{1}{t}F(tx,ty)$.

Suite aux travaux de Dixmier-Kirillov-Duflo on apprend que l'analyse sur les groupes et algèbres de
Lie ne porte pas sur les fonctions mais sur les demi-densités. On sait aussi depuis
Harish-Chandra que
l'objectif est  de ramener par l'application exponentielle l'analyse sur les groupes à
l'analyse sur les algèbres de Lie. Dans ce contexte la racine carrée du jacobien de l'application
exponentielle intervient inévitablement. Notons comme à l'habitude par $j(x)$ la fonction
$\det_{\g}(\frac{(1-\exp(-\mathrm{ad}x)}{\mathrm{ad}x})$ et par $q(x)$
la fonction
$\det_{\g}(\frac{\sinh(\mathrm{ad}x/2)}{\mathrm{ad}x/2})$. L'écart en ces  deux fonctions
n'est pas significatif dans cet article.

Il est bien connu  que l'on a la formule suivante \cite{KV}
\begin{equation}\label{eqdiffj}
j^{-1/2}(tx)\frac{\partial}{\partial t}j^{1/2}(tx)=\frac{1}{2}\mathrm{tr}_{\g}(\frac{\mathrm{ad}x}
{\exp(t\mathrm{ad}x)-1}-\frac{1}{t}).
\end{equation}

On sait aussi depuis longtemps qu'il ne faut pas considérer les distributions comme un module
à gauche sur les opérateurs différentiels mais comme un module à droite. Remarquons que ceci est
bien connu des analystes algébriques.
Un des buts de l'article de Kashiwara-Vergne était de démontrer que le produit
de convolution se transporte par l'application exponentielle quand on se restreint à des
distributions invariantes. Cette conjecture est démontrée dans sa totalité dans l'article
\cite{AST} (qui fait suite à \cite{ADS}) en utilisant les méthodes de Kontsevich. On peut considérer cet article comme une suite
naturelle de \cite{AST}. Plus précisément il s'agissait de montrer que si $u$ et $v$ sont
deux distributions invariantes au voisinage de $0$ dans $\g$  vérifiant une certaine condition
de support afin d'assurer un sens à la convolution, on avait la formule suivante

\begin{equation}\label{eqKV}
\int u(x)v(y)\frac{j^{1/2}(x)j^{1/2}(y)}{j^{1/2}(Z(x,y))}f(Z(x,y))dxdy=\int u(x)v(y)f(x+y)dxdy
\end{equation}
pour $f$ une fonction $C^{\infty}$ dans un voisinage de $0$. Cette
formule est démontrée dans \cite{AST}. La stratégie initiée dans
Kashiwara-Vergne était de  tenir compte de la déformation en le
paramètre $t$ et de démontrer que l'on avait la formule pour tout
$t$

\begin{equation}\label{eqKVt}
\int u(x)v(y)\frac{j^{1/2}(tx)j^{1/2}(ty)}{j^{1/2}(tZ_{t}(x,y))}f(Z_{t}(x,y))dxdy=\int u(x)v(y)f(x+y)dxdy.
\end{equation}
L'idée est maintenant simple, il suffit de dire que la dépendance en $t$
est triviale, c'est à dire que  la dérivée par rapport à $t$ est nulle.

Calculons cette dérivée comme dans \cite{KV}. Notons $D(tx,ty)$ la fonction de densité
$$\frac{j^{1/2}(tx)j^{1/2}(ty)}{j^{1/2}(tZ_{t}(x,y))}.$$
 On peut  remplacer $j$ par $q$ dans cette formule, la fonction de densité reste la même.
On a facilement compte tenu des équations (\ref{eqdiff}) et  (\ref{eqdiffj})

$$\frac{\partial}{\partial t}D(tx,ty)=\frac{1}{2}\mathrm{tr}_{\g}(\frac{\mathrm{ad}x}
{\exp(t\mathrm{ad}x)-1}+\frac{\mathrm{ad}y}{\exp(t\mathrm{ad}y)-1}
-\frac{\mathrm{ad}Z_t(x,y)}{\exp(t\mathrm{ad}Z_{t}(x,y))-1}-\frac{1}{t})D(tx,ty)+$$

\begin{equation}\label{deriveeA}
[x, F_{t}(x,y))]\cdot\partial_{x}D(tx,ty) +[y, G_{t}(x,y))]\cdot\partial_{y}D(tx,ty).
\end{equation}
Pour simplifier on va noter
\begin{equation}\label{T(x,y)}
T(x,y)=\frac{1}{2}\mathrm{tr}_{\g}(\frac{\mathrm{ad}x}
{\exp(\mathrm{ad}x)-1}+\frac{\mathrm{ad}y}{\exp(\mathrm{ad}y)-1}
-\frac{\mathrm{ad}Z(x,y)}{\exp(\mathrm{ad}Z(x,y))-1}-1).
\end{equation}
Par conséquent le terme qui apparaît en (\ref{deriveeA}) est
$\frac{1}{t}T(tx,ty)$. Ce calcul se justifie comme suit. Le
premier terme résulte de la dérivée par rapport à $t$ dans les
termes  $j^{1/2}(t\cdot)$ et le second terme résulte de la dérivée
en  $t$ dans  $Z_{t}$. Plus précisément compte tenu de
(\ref{eqdiff}) il vient que pour toute fonction $\phi$ on a
\begin{equation}
\frac{\partial}{\partial t} \phi(Z_{t}(x,y))=[x, F_{t}(x,y))]\cdot\partial_{x}\phi(Z_{t}(x,y)) +
[y, G_{t}(x,y))]\cdot\partial_{y}\phi(Z_{t}(x,y)
\end{equation}
Le  champ de vecteur $[x, F_{t}(x,y))]\cdot\partial_{x} +[y, F_{t}(x,y))]\cdot\partial_{y}$
agit trivialement sur la fonction
\\
$j^{1/2}(tx)j^{1/2}(ty)$ car cette dernière est invariante en chaque variable
sous l'action adjointe par conséquent la dérivée du terme en  $Z_{t}$ s'écrit bien comme annoncée.

On peut maintenant terminer le calcul de la dérivée dans (\ref{eqKVt}). Il vient

$$\frac{\partial}{\partial t}(D(tx,ty) f(Z_{t}(x,y)))=$$
$$[x, F_{t}(x,y))]\cdot\partial_{x}(D(tx,ty) f(Z_{t}(x,y))) +[y, G_{t}(x,y))]\cdot\partial_{y}(D(tx,ty) f(Z_{t}(x,y))) +$$
$$ \frac{1}{t}T(tx,ty)D(tx,ty) f(Z_{t}(x,y)).$$
On est donc amené  à calculer l'action à droite sur la distribution $u(x)v(y)$ du champ de vecteur
$[x, F_{t}(x,y))]\cdot\partial_{x} +[y, F_{t}(x,y))]\cdot\partial_{y}$. Compte tenu de l'invariance
de cette distribution on a
\begin{equation}\label{invarianceUV}
u(x)v(y)([x, F_{t}(x,y))]\cdot\partial_{x} +[y, F_{t}(x,y))]\cdot\partial_{y})=
-u(x)v(y)\mathrm{tr}_{\g}(\mathrm{ad}x\partial_{x}F_{t}(x,y)+\mathrm{ad}y\partial_{y}G_{t}(x,y)).
\end{equation}
Pour conclure au transport de la convolution dans  (\ref{eqKV})
 il suffit de demander que l'on ait
\begin{equation}
\frac{1}{t}T(tx,ty)-\mathrm{tr}_{\g}(\mathrm{ad}x\partial_{x}F_{t}(x,y)+\mathrm{ad}y\partial_{y}G_{t}(x,y))=0.
\end{equation}

La conjecture combinatoire de Kashiwara-Vergne est précisément cette égalité.

\subsection{Déformation de Kontsevich et résultats}
Dans son article fondamental  M. Kontsevich \cite{kont} introduit pour les variétés de Poisson une
déformation
formelle. Appliquée au cas du dual d'une algèbre de Lie on retrouve comme déformation de l'algèbre
symétrique le star produit de Duflo, à savoir
la symétrisation modifiée par l'opérateur $\partial(q^{1/2})$. D'une certaine manière on peut
penser que l'on a rien de neuf. Toutefois le point
intéressant est que les c\oe fficients apparaissant dans la formule de Campbell-Hausdorff se
calculent  à l'aide d'intégrales explicites. Compte tenu de la flexibilité naturelle de ces
intégrales et de l'argument d'homotopie expliqué dans Kontsevich et repris dans \cite{ADS},
\cite{AST} et \cite{MT}, on peut penser qu'il existe une version ``à la Kontsevich''
de la conjecture combinatoire.

Expliquons notre stratégie dans cette introduction.

En examinant la démonstration de Kashiwara-Vergne on s'aperçoit que le paramètre de
déformation crée  une dépendance en $t$ de la formule de Campbell-Hausdorff notée $Z_{t}$.
L'équation (\ref{eqdiff}) relie la différentielle en
$t$ de $Z_{t}$ à l'action d'un champ de vecteur.
La formule (\ref{deriveeA}) donne la dérivée en $t$ de la fonction de densité
 $D_{t}(x,y)=D(tx,ty)$.

D'une certaine manière on pourrait dire que la conjecture de Kashiwara-Vergne   est faite
pour  que l'on ait l'équation
\begin{equation}\label{eqA}
\frac{\partial}{\partial t}D_{t}(x,y)=
\mathrm{tr}_{\g}(\mathrm{ad}x\partial_{x}F_{t}+\mathrm{ad}y\partial_{y}G_{t})D_{t}(x,y)+[x, F_{t}(x,y))]\cdot\partial_{x}D_{t}(x,y) +[y, G_{t}(x,y))]\cdot\partial_{y}D_{t}(x,y)
\end{equation}
puisque l'autre équation
\begin{equation}
\frac{\partial}{\partial t} f(Z_{t}(x,y))=
[x, F_{t}(x,y))]\cdot\partial_{x} f(Z_{t}(x,y)) +[y, G_{t}(x,y))]\cdot\partial_{y} f(Z_{t}(x,y))
\end{equation}
résulte automatiquement de l'équation
\begin{equation}\label{eqZ}
\frac{\partial}{\partial t} Z_{t}(x,y)=
[x, F_{t}(x,y))]\cdot\partial_{x} Z_{t}(x,y) +[y, G_{t}(x,y))]\cdot\partial_{y} Z_{t}(x,y).
\end{equation}
Le problème cousin de Kashiwara-Vergne est donc celui-ci.
Existe-t-il une autre déformation de la formule de
Campbell-Hausdorff telle que les équations (\ref{eqZ}) et
(\ref{eqA}) soient vérifiées? On remarquera alors que ce résultat
entraîne aussi que l'exponentielle transporte la convolution des
distributions invariantes. En effet des équations (\ref{eqZ}) et
(\ref{eqA}) on déduit que l'on a

$$\frac{\partial}{\partial t}(D_{t}(x,y) f(Z_{t}(x,y)))=$$
$$[x, F_{t}(x,y))]\cdot\partial_{x}(D_{t}(x,y) f(Z_{t}(x,y))) +[y, G_{t}(x,y))]\cdot\partial_{y}(D_{t}(x,y) f(Z_{t}(x,y))) +$$
\begin{equation}
\mathrm{tr}_{\g}(\mathrm{ad}x\partial_{x}F_{t}+\mathrm{ad}y\partial_{y}G_{t})D_{t}(x,y) f(Z_{t}(x,y)).
\end{equation}
Comme précédemment, appliquée à des distributions invariantes $u(x)v(y)$ l'équation (\ref{invarianceUV}) assurera que le membre de gauche dans (\ref{eqKVt}) sera indépendant de $t$.

On peut maintenant énoncer notre théorème. On note $\overline{C}_{2,0}$ la variété de configuration de
deux points aériens dans le demi-plan de Poincaré (on reviendra sur cette variété introduite
par Kontsevich dans les  paragraphes suivants).
\\

\noindent {\bf{Théorème: }}  {\it Il existe une déformation $Z_{\xi}(x,y)$ de la formule de Campbell-Hausdorff $Z(x,y)$ dépendant d'un paramètre $\xi \in \ov{C}_{2,0}$, une déformation $D_{\xi}(x,y)$ de la fonction
de densité $D(x,y)$.  Il existe  des $1$-formes
différentielles en $\xi$ dépendant de  manière analytique de $(x,y)$ (dans un voisinage de $0$)
 à valeurs dans $\g$, on les note
 $\Omega_{F}(\xi,x,y)$ et
$\Omega_{G}(\xi,x,y)$. Les égalités suivantes sont alors vérifiées

$$dZ_{\xi}(x,y)=[x, \Omega_{F}(\xi, x,y)]\cdot\partial_{x} Z_{\xi}(x,y) +
[y, \Omega_{G}(\xi, x,y)]\cdot\partial_{y} Z_{\xi}(x,y)$$ et
$$dD_{\xi}(x,y)= \mathrm{tr}_{\g}(\mathrm{ad}x\partial_{x}\Omega_{F}(\xi,x,y)+
\mathrm{ad}y\partial_{y}\Omega_{G}(\xi,x,y)) D_{\xi}(x,y)$$
$$+[x, \Omega_{F}(\xi,x,y)]\cdot\partial_{x}D_{\xi}(x,y) +[y, \Omega_{G}
(\xi,x,y)]
\cdot\partial_{y}D_{\xi}(x,y).$$}

En appliquant ce théorème le long d'un chemin reliant le coin de l'oeil
 à l'iris (comme dans
\cite{ADS}) on retrouve une déformation dépendant d'un paramètre $t$
 comme dans
Kashiwara-Vergne. La conclusion est donc qu'il existe de nombreuses
déformations de la formule
de Campbell-Hausdorff qui fournissent automatiquement la condition des
traces. Malheureusement
ces déformations ne se font pas dans la catégorie des algèbres de Lie,
c'est à dire que la
fonction $Z_{\xi}(x,y)$ n'est pas la formule de Campbell-Hausdorff
d'une algèbre de Lie. En
d'autres termes on n'a pas
$$Z_{\xi}(Z_{\xi}(x,y),w)=Z_{\xi}(x,Z_{\xi}(y,w)).$$
Comme corollaire de notre théorème on obtient bien évidemment que
l'exponentielle transporte la convolution des distributions
invariantes (voir \cite{KV}, \cite{ADS} et \cite{AST} pour les
conséquences en algèbre et en analyse de ce théorème).

\section{Rappels sur la quantification de Kontsevich}

Dans cette section on va rappeler rapidement la construction de M.
Kontsevich pour la quantification formelle des variétés de Poisson
lorsqu'on l'applique au cas des structures de Poisson linéaires,
c'est à dire  dans le cas du dual des algèbres de Lie.
\subsection{Espaces de configurations}

On note par $C_{n,m}$ l'espace des configurations de $n$ points distincts
dans le
demi-plan de Poincaré (ce sont les points aériens) et $m$ points  distincts sur la droite réelle
(ce sont les points terrestres), modulo l'action
du groupe $az+b$ (pour $a \in \R^{+*}, b\in \R$). Compte tenu de l'action
de ce groupe sur les points réels, on peut identifier deux des points
réels aux points $0$ et $1$ (à condition que l'on ait $m\geq 2$, sinon on peut identifier
un des points aériens au complexe $i$). Dans son article \cite{kont} M. Kontsevich construit des
compactifications de ces variétés notées $\overline{C}_{n,m}$. Ce
sont des variétés à coins de dimension $2n+m-2$. Ces variétés ne sont pas connexes pour $m\geq 2$.
On notera par  $\overline{C}^{+}_{n,m}$ la composante  qui contient les configurations où les
points terrestres sont ordonnés dans l'ordre croissant (ie on  a $\ov{1}< \ov{2}<\cdots<\ov{m}$).

\begin{figure}[h!]
\begin{center}
\includegraphics[width=5cm]{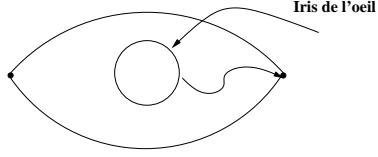}
\caption{\footnotesize L'oeil de Kontsevich $\overline{C}^{+}_{2,0}$ et son chemin}\label{oeil}
\end{center}
\end{figure}

\subsection{Graphes admissibles}
La notion de graphes  admissibles  est maintenant bien établie
dans la littérature. En l'occurrence on notera $G_{n,2}$ les
graphes avec $n$ points de première espèce et numérotés $1,2
\cdots n$ (points aériens) et $2$ points de seconde espèce (points
terrestres) numérotés $\overline{1}, \overline{2}$. Les graphes
qui interviennent dans le cas variétés de Poisson, sont des
graphes tels que de tout point aérien $i$ partent deux arêtes
numérotées $(e_{i}^{a},e_{i}^{b})$, mais il n'y a pas d'arêtes qui
bouclent et il n'y a pas d'arêtes doubles. On note $E_{\Gamma}$
l'ensemble des arêtes du graphe $\Gamma$.

\begin{figure}[!h]
\begin{center}
\includegraphics[width=5cm]{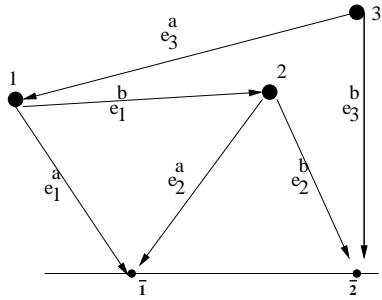}
\caption{\footnotesize Graphe admissible de type Lie}\label{graphe admissible}
\end{center}
\end{figure}

Dans le cas des algèbres de Lie, chaque sommet de première espèce
ne peut recevoir qu'au plus une arête: ce sont les graphes
pertinents (relevant graphs)  et on renvoie le lecteur à l'article
\cite{AST} pour une description précise des graphes pertinents
qui interviennent dans la formule finale.

Quand on a un produit de Poisson disons $\alpha$ sur $\R^{d}$, on peut alors associer à tout graphe admissible un
opérateur bi-différentiel sur $\R^{d}$ comme expliqué dans \cite{kont}.
On notera $B_{\G}(f,g)$ l'opérateur
bidifférentiel associé que l'on suppose agir sur les fonctions $f$ et $ g$. Expliquons par
souci de concision la formule. Sur chaque sommet aérien on met le crochet
de Poisson et sur les sommets terrestres on met les fonctions $f$ et $g$. Chaque arête arrivant sur un sommet dérive la fonction associée au sommet. On
multiplie les fonctions ainsi obtenues et on somme sur toutes les possibilités.
Concrètement la formule est la suivante. Pour chaque arête $e$, on note par $s(e)$ le point aérien source (départ) et par $b(e)$ le point but (arrivée).
Dans la formule ci-dessous $I$ décrit l'ensemble  les applications de l'ensemble des
arêtes $E_{\G}$ dans l'ensemble des indices de coordonnées $\{1, \cdots d\}$.

$$
B_{\Gamma,\alpha}(f,g) = \sum_{I} \bigg[\prod_{k=1}^n
\big(\prod_{e \in E_{\Gamma}, b(e) = k} \partial_{I(e)}\big)
\alpha^{I(e_k^a)I(e_k^b)}\bigg]\big( \prod_{e \in E_{\Gamma}, b(e) =
\ov{1}}\partial_{I(e)}\big) f \big(\prod_{e\in E_{\Gamma}, b(e) = \ov{2}}
\partial_{I(e)}\big) g.$$

\subsection{Forme d'angles}

Soient deux points distincts $(p,q)$ dans le demi-plan de Poincaré muni
de la métrique de Lobachevsky. On note
\begin{equation}
\phi(p,q)=\frac{1}{2i}\log(\frac{(q-p)(\overline{q}-p)}
{(q-\overline{p})(\overline{q}-\overline{p})}).
\end{equation}
C'est l'angle entre la géodésique $(p,\infty)$ et $(p,q)$ où l'infini
peut être vu comme l'infini sur la droite réelle. Quand on prend comme
modèle du demi-plan de Poincaré le disque unité alors l'infini en question est bien-sûr
 le point manquant sur le cercle unité. Géométriquement
la géodésique $(p,\infty)$ est tout simplement la demi-droite verticale
issue de $p$.

\begin{figure}[!h]
\begin{center}
\includegraphics[]{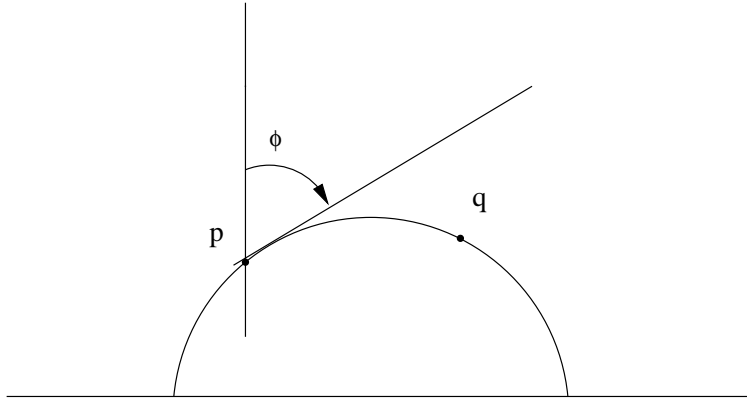}
\caption{\footnotesize Fonction d'angle}
\end{center}
\end{figure}

La fonction  d'angle s'étend en une fonction régulière à la compactification
 $\ov{C}_{2,0}$. Cette dernière est précisément décrite dans l'article de Kontsevich,
c'est le fameux ``oeil''(figure [\ref{oeil}]). On remarquera, mais c'est tautologique vu la
construction des compactifications que lorsque les points $p, q$ s'approchent selon
un angle $\theta$, la fonction d'angle vaut précisément cet angle. Lorsque $p$ s'approche
de l'axe réel, la fonction d'angle est nulle et lorsque c'est $q$ qui s'approche de l'axe
réel on obtient deux fois l'angle de demi-droite avec l'axe réel.

\begin{figure}[h!]
\begin{center}
\includegraphics[height=6cm, width=8cm]{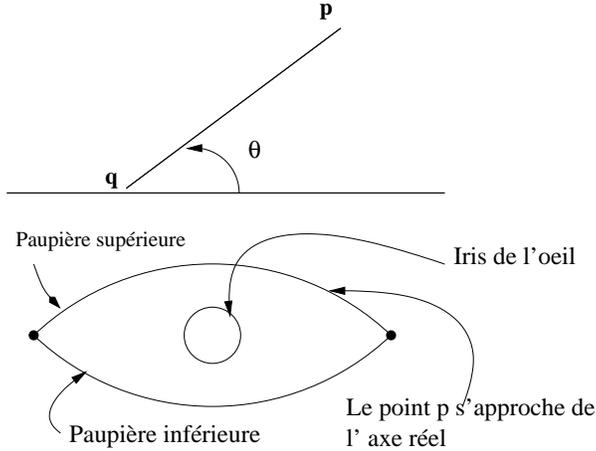}
\caption{\footnotesize Valeur de la fonction d'angle le long de la paupière}
\end{center}
\end{figure}

Dans \cite{kont} il est expliqué que  l'on peut prendre n'importe qu'elle fonction
d'angle pour peu qu'elle vérifie certaines propriétés, mieux encore on peut
même préférer  une autre forme différentielle \cite{operades}.

Bref, comme la fonction est régulière sur la compactification, on peut considérer
sa différentielle qui est alors une $1$-forme sur $\overline{C}_{2,0}$.

\subsection{Poids associé à un graphe}
Si $\Gamma$ est un graphe dans $G_{n,2}$, alors  toute arête $e$ définit par restriction
une fonction d'angle notée $\phi_{e}$ sur la variété $\overline{C}^{+}_{n,2}$. Le produit ordonné
\begin{equation}
\Omega_{\Gamma}=\bigwedge _{e \in E_{\Gamma}} d\phi_{e}
\end{equation}
est donc une $2n$-forme sur  $\overline{C}^{+}_{n,2}$ qui est de dimension aussi $2n$.
Le poids associé sera donc

\begin{equation}
w_{\Gamma}=\frac{1}{(2\pi)^{2n}}\int_{\overline{C}^{+}_{n,2}} \Omega_{\Gamma}.
\end{equation}

\subsection{Formule de Quantification}
Dans le cas  des structures de Poisson sur
$\R^{d}$, la formule de Kontsevich
s'écrit en terme des ingrédients introduits de la façon suivante.
Pour $f$ et $g$ deux fonctions polynomiales sur $\R^{d}$ on a
\begin{equation}\label{formuleK}
f\star_{h} g=fg +\sum _{n=1}^{n=\infty}\frac{h^{n}}{n!}\sum_{\G \in G_{n,2}}w_{\Gamma} B_{\G}(f,g)
\end{equation}
avec $h$ un paramètre formel. Dans le cas des algèbres de Lie  l'ordre des opérateurs
$B_{\G}$ est suffisamment croissant et  la formule ci-dessus est en fait une somme finie.
En faisant $h=1$ ce star produit vérifie pour $f=x$ et $g=y$ la relation souhaitée
$$x\star y-y\star x=[x,y]$$
lorsque l'on prend pour crochet de Poisson la moitié du crochet de Lie.

\section{La formule Campbell-Hausdorff}
\subsection{La formule Campbell-Hausdorff en termes de diagrammes}

\begin{figure}[h!]
\begin{center}
\includegraphics[width=5cm]{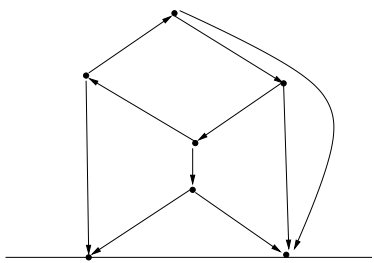}
\caption{\footnotesize Graphe simple de type Roue}\label{roue}
\end{center}
\end{figure}

Comme dans \cite{AST} ou \cite{ka} la formule de Kontsevich
appliquée au cas des algèbres de Lie donne le résultat suivant.
Tout graphe admissible pertinent se décompose en produit de
graphes simples. Chaque graphe simple est soit de type Lie  i.e.
un graphe avec une seule racine et sans roue (le graphe est un
arbre comme dans la figure [\ref{graphe admissible}]) soit un
graphe de type roue (ie une roue tentaculaire comme dans la figure
[\ref{roue}]).

Dans  \cite{AST} on a associé un symbole à chaque graphe, on
notera abusivement $\Gamma(x,y)$ le symbole associé au graphe
$\Gamma$. Par exemple le symbole associé au graphe de la figure[
\ref{graphe admissible}] est $\G(x,y)=[[x,[x,y]],y]$ et le symbole
associé au graphe de la figure [\ref{roue}] est $\tr( \ad x
\ad[x,y]\ad y\ad y)$. Le symbole est  une fonction polynomiale de
$\g\times\g$ dans $S[\g]$.

Lorsque le graphe est simple et de type Lie, alors $\Gamma(x,y)$
est naturellement un élément de l'algèbre de Lie engendrée par $x$
et $y$.

Comme démontré dans \cite{ka}, la formule Campbell-Hausdorff s'écrit pour $x$ et $y$ dans $\g$

\begin{equation}
Z(x,y)=x+y +\sum_{\Gamma}w_{\Gamma}\Gamma(x,y)
\end{equation}
où la somme porte sur l'ensemble des graphes ``géométriques''
simples et de type Lie. Les graphes qui contribuent de manière non
triviale dans cette formule n'ont donc qu'une racine et ne
possèdent pas de symétries. Par conséquent les graphes numérotés
associés à un graphe géométrique sont au nombre de $n!2^n$. Le
terme $n!$ disparaît dans la formule finale car dans la formule de
Kontsevich (\ref{formuleK})
 on avait un terme en $1/n!$ et le terme
$2^n$ disparaît aussi car on prend dans le cas des algèbres de Lie
pour crochet de Poisson $\frac{1}{2}$ du crochet de Lie. Toutefois
le symbole $\Gamma(x,y)$ est mal défini si  le graphe n'est pas
numéroté. Pour résoudre ce problème il suffit de remarquer que
c'est aussi le cas pour le c\oe fficient $w_{\Gamma}$ et  que les
deux difficultés se compensent.

\subsection{Commentaires}

Les symboles $\Gamma(x,y)$ ne forment
pas une base de l'algèbre de Lie libre engendrée par $x$ et $y$, par conséquent il y a des
redondances dans la formule, mais c'est aussi le cas quand on utilise d'autres formules donnant
Campbell-Hausdorff.

La nature des c\oe fficients $w_{\Gamma}$ reste mystérieuse. Les
termes associés aux produits du genre $[x,[x,[x,\cdots,[x,y]]]]$
sont des nombres de Bernoulli comme tout le monde le sait, puisque
ces termes calculent la différentielle de l'application
exponentielle. C'est aussi le cas pour les c\oe fficients de
Kontsevich comme remarqué dans \cite{ka}. La différentielle de
l'application exponentielle vaut
$$\frac{1-\exp(-adx)}{adx},$$
il vient que l'on a à l'ordre $1$ en $y$
\begin{equation}
Z(x,y)=x+ \frac{adx}{1-\exp(-adx)}\cdot y.
\end{equation}
Par ailleurs il n'y a qu'une seule façon d'écrire un terme de type
Lie où n'intervient qu'une seule fois $y$. Ce point est essentiel,
car il permet de conclure que les c\oe fficients $w_{\Gamma}$ sont
aussi de Bernoulli. Ce n'est plus le cas pour les autres. On ne
peut pas a priori affirmer que les nombres $w_{\Gamma}$ sont tous
rationnels  et je ne sais pas si c'est le cas. Quand les graphes
sont plus généraux on sait que ces nombres ne sont pas rationnels
\cite{operades}. La conclusion est que
 les termes non rationnels se compensent à cause de l'identité de Jacobi.
En conséquence on considère une écriture de Campbell-Hausdorff qui
 a le gros inconvénient de n'être pas
a priori une expression rationnelle (bien-sûr la somme peut se réécrire sous forme rationnelle!).

Compte tenu de cette remarque, on voit que la déformation que l'on va proposer est largement
différente de celle proposée dans Kashiwara-Vergne, qui restait rationnelle et
dans la catégorie des algèbres de Lie.

\subsection{Analyse des termes d'ordre $1$ dans l'équation de Kashiwara-Vergne}

Compte tenu de l'équation
\begin{equation}
\frac{\partial}{\partial t}Z_{t}(x,y)=[x, F_{t}(x,y))]\cdot\partial_{x}Z_{t} +[y, G_{t}(x,y))]\cdot\partial_{y}Z_{t},
\end{equation}
 d'où viennent les termes d'ordre $1$ en $y$ dans la dérivée
$\frac{\partial}{\partial t}Z_{t}(x,y)$?

Par symétrie Kashiwara-Vergne montrent  que l'on peut prendre

$F(x,y)=G(-y,-x)$.
On voit facilement que l'on a
$F(x,y)=\frac{1}{4}y+\cdots$ et $G(x,y)=-\frac{1}{4}x+\cdots$.

Ecrivons pour simplifier à l'ordre $1$ en $y$

\begin{equation}
Z(x,y)=x+ \frac{adx}{1-\exp(-adx)}\cdot y=x+\sum_{n\geq 0} b_n \ad^{n}(x)y.
\end{equation}
il vient

\begin{equation}
\frac{\partial}{\partial t}Z_{t}(x,y)=\sum_{n\geq  1} n t^{n-1} b_n \ad^{n}(x)y
\end{equation}
Or le terme $[x, F_{t}(x,y))]\cdot\partial_{x}Z_{t}(x,y)$ contribue
(à l'ordre $1$ en $y$ toujours) comme
 $[x, F_{t}(x,y))]\cdot\partial_{x} x$ c'est à dire $[x, F_{t}(x,y))]$ (car le terme en $y$ est forcément dans le crochet $[x, F_{t}(x,y))]$). Ecrivons
$$F(x,y)=\sum d_{n}t^{n}\ad^{n}(x)y$$
avec $d_{0}=\frac{1}{4}$.

Le terme $[y, G_{t}(x,y))]\cdot\partial_{y}Z_{t}(x,y)$ contribue comme
$$[y, -\frac{1}{4}x]\cdot\partial_{y}\sum b_n t^{n}\ad^{n}(x)y=
\sum \frac{b_{n}t^{n}}{4} \ad^{n}(x)[x,y].$$ On en déduit la
relation de récurrence
$$(n+1)b_{n+1}=d_{n} +  \frac{b_n}{4} $$
avec $$b_{0}=1, b_{1}=\frac{1}{2}, b_{2}=\frac{1}{12}, b_{3}=0, b_{4}=\frac{1}{720}.$$
Il vient alors
$$d_{0}=\frac{1}{4},d_{1}=\frac{1}{24},d_{3}=-\frac{1}{48},d_{4}=-\frac{1}{180}.$$

La situation sera comme on va le voir beaucoup plus directe dans notre formule mais restera quand
même dans cet esprit. Remarquons aussi que l'on a pu faire le calcul ci-dessus car il y a
unicité dans l'écriture des éléments de type Lie ne contenant qu'une seule fois $y$.

\section{Déformation des c\oe fficients dans Campbell-Hausdorff}

L'idée est très simple. Au lieu de déformer le crochet de Lie comme dans
Kashiwara-Vergne, on va déformer les c\oe fficients. Ceci est possible et
naturel dans la mesure où ces derniers se calculent grace à des intégrales.

Plus précisément, au lieu de prendre un graphe $\G$ dans $G_{n,2}$ prenons
ce même graphe géométrique mais dans $G_{n+2,0}$. C'est l'argument de base de
la déformation homotopique dans Kontsevich. On fixe donc un point $\xi$ dans
$\overline{C}^{+}_{2,0}$ représentant les points $n+1, n+2$. On peut alors
intégrer la forme différentielle associée au nouveau graphe lorsqu'on suppose
que le point $\xi$ est fixe. En terme simple, on intègre selon les $n$ premiers points
ce qui fait toujours $2n$ dimensions. On pose alors $w_{\G}(\xi)$ pour
le résultat obtenu. Lorsque le point $\xi$ vaut le point $(0,1)$ on retrouve
le c\oe fficient de Kontsevich. La dépendance en $\xi$ est régulière au sens des variétés
 à coins et on dispose d'une majoration uniforme en $\xi$ de ce c\oe fficient
comme dans \cite{ADS}.

\begin{figure}[h!]
\begin{center}
\includegraphics[height=4cm]{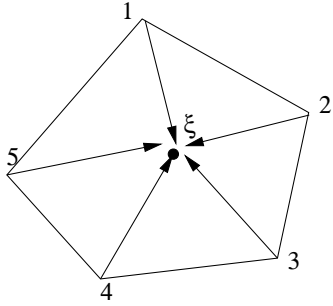}
\caption{\footnotesize Graphe associé à une roue pure dans le demi-plan de Poincaré}\label{roue de velo}
\end{center}
\end{figure}

Lorsque le point $\xi$ est sur l'iris de l'oeil c'est à dire que
les deux points sont infiniment proche selon un angle donné, on
trouve $0$. En effet les seuls c\oe fficients non nuls seraient
ceux correspondant à des roues de vélo (figure [\ref{roue de
velo}]) avec les rayons dirigées vers le point $\xi$, mais
Shoikhet \cite{shoi} a montré que ces termes sont aussi nuls dans
le cas des algèbres de Lie.

On dispose maintenant comme dans \cite{ADS} d'une  produit non associatif
\begin{equation}
f\star_{\xi}g=\sum _{n, \G\in G_{n,2}}\frac{w_{\G}(\xi)}{n!}B_{\G}(f,g)
\end{equation}
où on place la fonctions $f$ en $n+1$ et la fonction $g$ en $n+2$.
La formule commence par le produit $fg$ si l'on admet que pour $n=0$ le poids associé au graphe sans point aérien  est $1$, mais il est plus commode d'écrire
\begin{equation}
f\star_{\xi}g=fg +\sum _{n\geq 1, \G\in G_{n,2}}\frac{w_{\G}(\xi)}{n!}B_{\G}(f,g).
\end{equation}

 Par suite on dispose d'une  déformation de la formule de Campbell-Hausdorff
\begin{equation}
Z_{\xi}(x,y)=x+y+\sum_{\G} w_{\G}(\xi)\G(x,y)
\end{equation}
où la somme porte sur les graphes simples de type Lie c'est à dire les graphes
géométriques avec une seule racine. Dans cette formule on n'a plus besoin
 de considérer des graphes numérotés. C'est ce que l'on a
appelé les graphes géométriques.
\\
Compte tenu des résultats de Shoikhet on a $Z_{\xi \in
\mathrm{Iris}}(x,y)=x+y$ et
 $Z_{\xi=(0,1)}(x,y)=Z(x,y)$.

Le point nouveau ici et c'est plus ou moins l'objet de cet article, c'est que ce nous allons calculer directement la
dérivée par rapport à $\xi$. C'est naturel et c'est de manière indirecte ce
qui fut déjà fait dans \cite{ADS} et \cite{AST}. On va voir qu'il vient
une équation différentielle.

\subsection{Calcul de la dérivée}

Lorsque l'on veut calculer la différentielle  de $w_{\xi}$ on écrit $w_{\nu}-w_{\xi}$
et on fait tendre $\nu$ vers $\xi$, ce qui revient à
prendre un chemin infinitésimal dans l'oeil de Kontsevich de $\xi$ vers $\nu$.
Le raisonnement de base fait de nombreuses fois est donc d'écrire la formule
de Stokes comme dans \cite{ADS} mais cette fois sur un chemin infinitésimal.

Les mêmes types de composantes sont à étudier. On ne doit regarder, pour des
raisons de dimension (voir \cite{MT} pour une étude plus précise sur le cup-produit), que les regroupements de deux points aériens. En bref le calcul de la dérivée des c\oe fficients
s'obtient en contractant les arêtes.

Toutefois lorsqu'on s'intéresse à la somme partielle
$$\sum_{\G } w_{\G}(\xi)\G(x,y)$$
où $\G$ est un graphe simple de type Lie avec $n+2$ sommets aériens dont $2$
représentent le point $\xi$, il y aura une simplification à cause de l'identité de Jacobi.

\begin{figure}[!h]
\begin{center}
\includegraphics[width=5cm]{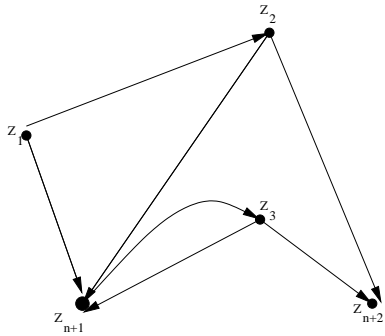}
\caption{\footnotesize Graphe type obtenu quand un point s'agglutine sur $z_{n+1}$}\label{graphe type}
\end{center}
\end{figure}

Décrivons les différents termes significatifs  qui apparaissent dans le calcul de la dérivée.

\begin{enumerate}
\item
Lorsque deux points autres que $(z_{n+1},z_{n+2})=\xi$ se
rapprochent,  l'identité de Jacobi montre que la contribution de
ce terme est nul. On remarquera ici que la nullité est dû à
l'identité de Jacobi pour les symboles $\G(x, y)$.
\item
Lorsque un point s'approche de $z_{n+1}$ on obtient un graphe type  comme à la figure
[\ref{graphe type}].

\item
Lorsque un point s'approche de $z_{n+2}$ on obtient une situation analogue.
\end{enumerate}

Remarquons encore une fois qu'il y a d'autres composantes de bord (des points aériens pourraient
s'agglutiner sur l'axe réel par exemple).
mais celles-ci ne contribuent pas dans la dérivée du c\oe fficient.

\begin{figure}[!h]
\begin{center}
\includegraphics[width=10cm]{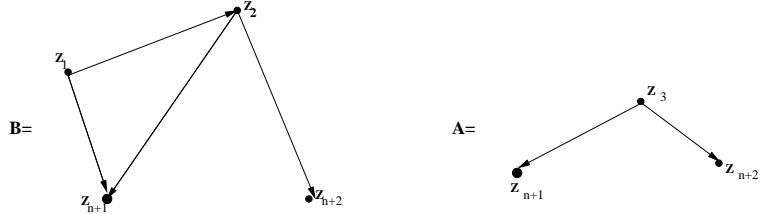}
\caption{\footnotesize Le graphe type de la figure [\ref{graphe type}] se décompose naturellement}\label{decompositionfig}
\end{center}
\end{figure}

Il faut comprendre maintenant que les arêtes issue de $z_{n+1}$ ou
arrivant sur  $z_{n+1}$ ou $z_{n+2}$ contribuent pour la
 forme différentielle totale comme une dérivée sur $\xi$. L'arête issue de
 $z_{n+1}$ va sur la racine d'un sous graphe de type Lie (et simple forcément).
On va noter $A$ ce sous-graphe et  on note $B$ le graphe obtenu en enlevant le sous-graphe
$A$ et l'arête qui joint  $z_{n+1}$ à $A$ (voir figure [\ref{decompositionfig}]).
La forme associée au graphe type
est une $(2n-2)$-forme à valeurs
dans les $1$-formes en $\xi$. On peut donc intégrer cette $(2n-2)$-forme sur les
configurations de $n-1$ points dans le demi-plan de Poincaré. On note $w_{A,B}^{x}(\xi)$
la $1$-forme obtenue après intégration. L'exposant $x$ est là pour rappeler qu'un point
s'est agglutiné sur $z_{n+1}$ et qu'il y  a une arête du point $z_{n+1}$ vers la racine de $A$.
De manière analogue on laissera un exposant $y$ dans
le cas où un point s'agglutine sur $z_{n+2}$. La proposition clé est  la suivante.
\begin{prop}\label{produit}
On a  $w_{A,B}^{x}(\xi)=w_{B}(\xi)w_{A}^{x}(\xi)$.
\end{prop}
{\it Preuve}:
Il faut remarquer d'abord que les graphes $A$ et $B$ sont indépendants, c'est
à dire que les seuls sommets d'intersection sont les points $z_{n+1}, z_{n+2}$.
Cela résulte de \cite{AST} sur la nature géométrique des graphes dans le cas
des algèbres de Lie.

Maintenant les termes qui contribuent pour la $1$-forme en $\xi$ proviennent
soit de $A$ soit de $B$. Considérons une contribution issue de $B$, on garde
alors une des arêtes arrivant sur $z_{n+1}$ ou $z_{n+2}$ comme différentielle en
$\xi$, ce qui montre qu'il va nous manquer une dimension pour intégrer sur les
configurations des sommets de $B$. La conclusion est que pour les configurations
 des points de $B$ on intègre à $\xi$ constant, ce qui donne $w_{B}(\xi)$. Ce
qui reste est alors le terme $w_{A}^{x}(\xi)$ et qui contient toute la partie
différentielle en $\xi$.$\Box$

Précisons ici que $w_{A}^{x}(\xi)$ est bien une 1-forme en $\xi$ car il y a une arête de $z_{n+1}$ vers la racine du graphe $A$ et que l'on intègre sur les configurations des sommets de $A$.

On  regroupe maintenant tous les graphes simples $\G$ qui dégénèrent de la même façon
et contribuent comme $w_{A,B}^{x}(\xi)$.

Comment sont-ils obtenus?

Tout simplement en déployant les arêtes de $B$ qui arrivent sur  $z_{n+1}$.
En effet si on déployait une arête de $A$ arrivant sur $z_{n+1}$, soit
 le graphe obtenu ne serait pas simple (puisqu'il aurait deux racines), soit
le graphe aurait une roue ce qui est exclu.

En examinant ce procédé au niveau des symboles, on constate que la
contribution dans la somme partielle $\sum
w_{A,B}^{x}(\xi)\G(x,y)$ est précisément
\begin{equation}\label{termeproduit}
w_{A,B}^{x}(\xi) [x,A(x,y)]\cdot \partial_{x}B(x,y)
\end{equation}

On a un terme analogue portant sur $y$.

Notons  $\Omega_{F}(x,y)$ la $1$-forme $\sum_{A}w_{A}^{x}(\xi)A(x,y)$ la somme portant sur les graphes géométriques de type Lie et de manière analogue  $\Omega_{G}(x,y)=\sum_{A}w_{A}^{y}(\xi)A(x,y)$
( ces séries n'ont donc pas de termes constants).

On peut maintenant énoncer notre premier théorème
\begin{theo}\label{deformationCBH}
Les $1$-formes  $\Omega_{F}$ et $\Omega_{G}$ en $\xi\in \overline{C}_{2,0}$
sont données par des séries (de type Lie) convergentes  à valeurs dans $\g$ et vérifient

$$dZ_{\xi}(x,y)=[x,\Omega_{F}(x,y)]\cdot \partial_{x}Z_{\xi}(x,y)+[y,\Omega_{G}(x,y)]\cdot \partial_{y}Z_{\xi}(x,y)$$
\end{theo}
{\it Preuve}: Le fait que les séries soient convergentes en $x$ et $y$
dans un voisinage de $0$ résulte comme dans
 \cite{ADS} ou \cite{AST} de  l'estimation des c\oe fficients (les formes à intégrer sont
 régulières sur des variétés à bord compactes).
Compte tenu de la proposition précédente on  a
$w_{A,B}^{x}(\xi)=w_{B}(\xi)w_{A}^{x}(\xi)$, vient alors

$$\sum_{A,B} w_{A,B}^{x}(\xi) [x,A(x,y)]\cdot \partial_{x}B(x,y)=
\sum_{A,B}w_{B}(\xi)w_{A}^{x}(\xi)[x,A(x,y)]\cdot \partial_{x}B(x,y)=$$
$$[x,\Omega_{F}]\cdot \partial_{x}Z_{\xi}(x,y).$$
Par symétrie on a le résultat souhaité.$\Box$

{\bf Remarque:} Si $\gamma: [0,1] \mapsto \ov{C_{2,0}}$ est un
chemin reliant le coin de l'oeil à l'iris (cf figure
[\ref{oeil}]), alors la fonction $t \mapsto Z_{\gamma(t)}(x,y)$
est une déformation de la formule de Campbell-Hausdorff. Les c\oe
fficients $\int_{\gamma} w_{A}^{x}$ vérifie les mêmes majorations
que dans  \cite{ADS} ce qui explique l'analycité en $(x,y)$ pour
les formes $\Omega_{F}(x,y)$ et $\Omega_{G}(x,y)$
\subsection{Commentaires et exemples}

\subsubsection{Graphe de Bernoulli}

Notre but est de faire $\xi=(0,1)$ dans l'équation différentielle obtenue dans le
théorème. On aura alors $Z_{\xi=(0,1)}(x,y)=Z(x,y)$.

Regardons ce que donne le calcul de la différentielle sur les c\oe
fficients qui déforment les nombres de Bernoulli apparaissant dans
le calcul à l'ordre $1$ en $y$. Rappelons que dans la déformation
de Kashiwara-Vergne on avait déformé ces nombres en les
multipliant par $t^n$.
\\

\begin{figure}[!h]
\begin{center}
\includegraphics[width=6cm]{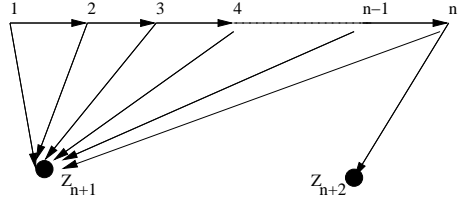}
\caption{\footnotesize Graphe type Bernoulli
déformé}\label{bernouilli}
\end{center}
\end{figure}

La déformation est donnée par le graphe de la figure [\ref{bernouilli}].
On note pour simplifier  $w_{n}(\xi)$ le c\oe fficient déformé.

\begin{figure}[!h]

\begin{center}
\includegraphics[width=6cm]{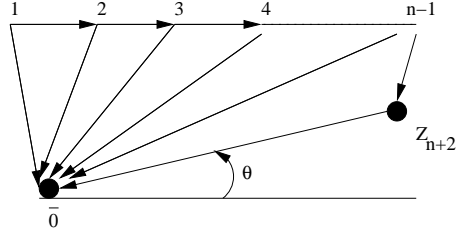}
\caption{\footnotesize Graphe calculant $w'_{n}(\theta)$}\label{derivee bernouilli}
\end{center}
\end{figure}

Le long de la paupière supérieure de l'oeil c'est à dire quand
 le point $z_{n+1}$ est sur l'axe réel, le c\oe fficient
  vérifie une équation différentielle très simple.
En effet on a alors $w_{n}(\xi)=w_{n}(\theta)$ avec $\theta$
l'angle entre $\ov{0}$ et $z_{n+2}$. La dérivée est donc
représentée par le graphe
 de la figure [\ref{derivee bernouilli}]. En effet toutes les autres contractions
 d'arêtes vont donner $0$ car le point $z_{n+1}$ reste sur l'axe réel. Seule
intervient la contraction de l'arête qui joint $z_{n}$ à $z_{n+2}$.

Or la partie différentielle en $\theta$ provient de l'arête joignant
 $\ov{0}$ et $z_{n+2}$ et qui représente clairement $\frac{d\theta}{\pi}$.
 Comme on l'a dit le reste des intégrations se fait à $\theta$ fixe, il vient donc

$$\frac{d}{d\theta}w_{n}(\theta)=-w_{n-1}(\theta)\frac{d\theta}{\pi}.$$
Ceci est conforme au résultat de \cite{ka} qui montre qu'on a
$$w_{n}(\theta)=\frac{(-1)^{n}}{n!}b_{n}(\frac{\theta}{\pi})$$ et il est bien connu
les polynômes de Bernoulli $b_{n}(\theta)$ vérifient
$$b_{n}'(\theta)=nb_{n-1}(\theta).$$

\subsubsection{Contribution d'ordre $1$ en $y$  au point $(0,1)$}

En cherchant les termes d'ordre $1$ en $y$ dans $\mathrm{d}
Z_{\xi}(x,y)$ on trouve alors deux types de contribution
représentée par les graphes suivants qui sont clairement les
seules façons de dégénérer les graphes associés au nombre de
Bernoulli.

\begin{figure}[!h]
\begin{center}
\includegraphics[width=6cm]{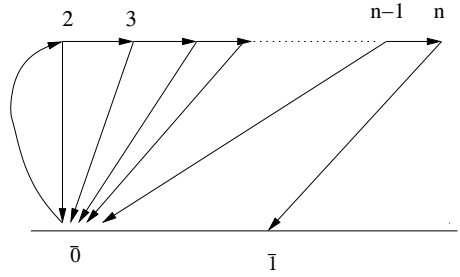}
\caption{\footnotesize Graphe calculant la dérivée au point $(0,1)$ le long de la paupière inférieure}\label{deriveeinferieure}
\end{center}
\end{figure}

\begin{figure}[!h]
\begin{center}
\includegraphics[width=6cm]{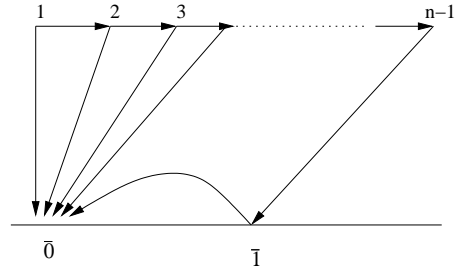}
\caption{\footnotesize Graphe calculant la dérivée au point $(0,1)$ le long de la paupière supérieure}\label{deriveesuperieure}
\end{center}
\end{figure}

Le premier (figure [\ref{deriveeinferieure}]) correspond d'après le théorème \ref{deformationCBH}
 et la formule (\ref{termeproduit}) à un terme du genre
$$[x,\ad (x)^{n-1}\cdot y]\partial_{x}(x),$$ c'est la dérivée le long
de la paupière inférieure prise au point $(0,1)$. On prendra garde que l'arête issue de $\ov{0}$
ne donne pas $0$ car le point $z_{n+1}$ s'approche de l'axe réel mais pas forcément le long
de l'axe réel.

Le deuxième (figure [\ref{deriveesuperieure}]) correspond à une contribution du genre
$$[y, x]\partial_{y}(\ad (x)^{n-1}\cdot y)$$
qui est la dérivée le long de la paupière supérieure prise au point $(0,1)$.

On peut dire la chose suivante. La déformation des c\oe fficients
est la moyenne des déformations le long des paupières. La
déformation le long de la paupière supérieure est tout simplement
la déformation des nombres de Bernoulli en les polynômes de
Bernoulli, l'autre déformation reste encore mystérieuse.

\begin{figure}[!h]
\begin{center}
\includegraphics[width=4cm]{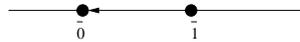}
\caption{\footnotesize Graphe de base}\label{graphe de base}
\end{center}
\end{figure}

\begin{figure}[!h]
\begin{center}
\includegraphics[width=4cm]{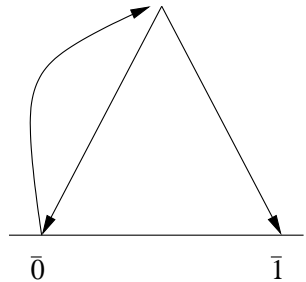}
\caption{\footnotesize Exemple de calcul de la dérivée}\label{graphe elem}
\end{center}
\end{figure}

Un calcul à la main, montre que si on normalise nos calculs
tel que le graphe de la figure [\ref{graphe de base}]
corresponde à $-1$, alors la dérivée correspondant
 au graphe de la figure [\ref{graphe elem}] vaut $1/4$.

Toutefois ce calcul semble montrer que la déformation
que nous proposons est différente au point $(0,1)$ de celle de Kashiwara-Vergne. En effet la
différentielle $\mathrm{d}Z_{\xi}$ en ce point n'a pas le même développement que
 $\dif{Z_{t}(x,y)}{t}$ pour $t=1$.

\subsubsection{Autre exemple  de calcul des c\oe fficients}

L'argument d'homotopie infinitésimal que nous avons utilisé pour
établir d'équation différentielle vérifiée par la formule
Campbell-Hausdorff, est clairement une application de la formule
de Stokes. Nous expliquons maintenant une procédure pour calculer
certains c\oe fficients associés à un graphe admissible $\G\in
G_{n,2}$. Nous renvoyons le lecteur à un prochain article dans
lequel nous développerons cette technique concernant les poids
associés à des graphes pertinents.

La procédure est la suivante. Ajoutons une arêtes au graphe. On a donc ajouté un point
aérien et une arête. Il y a un écart d'une dimension entre la variété d'intégration et
la forme différentielle associée au graphe. Appliquons la formule de
Stokes pour ce graphe $\tilde{\G}$. On a alors la formule
\begin{equation}
\int_{\ov{C}^{+}_{n+1,2}} \mathrm{d} \Omega_{\tilde{\G}}=
\int_{\partial \ov{C}^{+}_{n+1,2}}\Omega_{\tilde{\G}}.
\end{equation}

Cela revient à écrire  que la somme des c\oe fficients associés aux graphes  obtenus par
contraction d'une arête est nulle. On obtient alors en général une relation non triviale
entre c\oe fficients.

\begin{figure}[!h]
\begin{center}
\includegraphics[width=8cm]{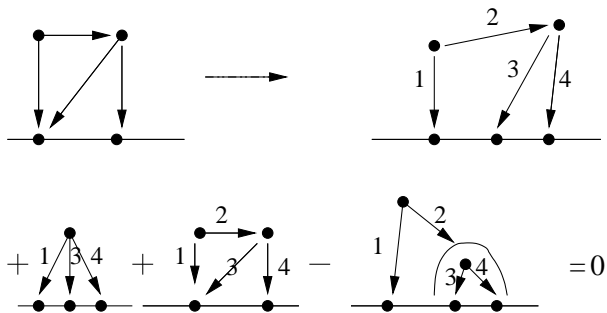}
\caption{\footnotesize Utilisation de Stokes dans le calcul des poids}\label{calculstokes2}
\end{center}
\end{figure}

Une autre méthode consisterait à ajouter un point terrestre  et à
déplacer certaines arêtes arrivant sur $\ov{0}$ par exemple  vers
ce point. On utilise alors la formule de Stokes, mais on prendra
garde qu'il faut dans ce cas bien étudier la contribution de
toutes les composantes de bord en particulier celles qui
correspondent aux regroupement de point aériens sur l'axe réel
(voir \cite{MT} pour plus de détails) et il faut suivre aussi les
orientations des faces (voir  \cite{AMM} pour des explications
détaillées).  Illustrons ceci par un exemple donné à la figure
[\ref{calculstokes2}].

On retrouve alors l'identité $1/6 +1/12 - (1/2)^2=0$, car le diagramme avec un point aérien et $3$
arêtes vaut $1/3!=1/6$.


\subsection{La déformation $Z_{\xi}$ et Campbell-Hausdorff}

Comme on l'a annoncé en introduction, notre déformation n'est pas associative en général, c'est
à dire que ce n'est pas une formule
de Campbell-Hausdorff d'une algèbre de Lie.
En effet une expression donnée par des polynômes de Lie qui
 définit une formule associative est forcément la formule de Campbell-Hausdorff
 de la structure de Lie
déterminée par le terme en $[x,y]$. C'est démontré dans
\cite{ABM}, mais le résultat est sans doute connu depuis longtemps
et on peut démontrer ceci par récurrence, en examinant les termes
de type Bernoulli. En particulier, il suffit de vérifier dans la
formule de $Z_{\xi}$,
 que les c\oe fficients associés aux termes de
type Bernoulli, ne vérifient pas les relations escomptées. Comme
le long de la paupière, ces c\oe efficients sont  les polynômes de
Bernoulli (au signe près), il suffit de constater que la relation
$w_{1}^{2}(\theta)=3w_{2}(\theta)$ n'est pas vérifiée. Par suite,
cette relation n'est pas vérifiée sur un ouvert  de la variété
$\ov{C_{2,0}}$.

\section{Déformation de la fonction de densité}
\subsection{Expression en terme de graphes}

La fonction de densité est la fonction $D(x,y)$ que nous avons
introduite au début de cet article. Comme démontré dans \cite{AST}
cette fonction s'exprime simplement en terme de symboles de
graphes. On a le résulta suivant
\begin{prop}
La fonction de densité $D(x,y)$ vérifie
$$A(x,y)=\sum_{\G\in G_{n,2}^{w}}\frac{1}{n!}w_{\G}\G(x,y)$$
où la somme porte sur tous les graphes admissibles de type roue, id est sans composante
simple de type Lie (par convention la série commence par 1).
\end{prop}
Notons que la somme portent sur tous les ``produits'' de graphes simples de type roue et non pas
seulement sur les graphes simples de type roue.
Comme on l'a remarqué dans \cite{AST} cela implique que cette fonction est l'exponentielle de
\begin{equation}\label{equationroue}
\sum_{\G\in G_{n,2}^{w}, \G \mathrm{simple}}w_{\G}\G(x,y)
 \end{equation}
mais où la somme porte uniquement sur les graphes simples géométriques de type
roue (figure [\ref{roue}]). Bien sûr cette série est convergente dans un voisinage de $0$ en $x$ et
 $y$.

\subsection{Déformation de la fonction de densité}

Effectuons la même procédure de déformation des c\oe fficients de type roue que dans le cas
Campbell-Hausdorff.  On note par analogie $D_{\xi}(x,y)$ la fonction
de densité obtenue par déformation des c\oe fficients.
 Il est plus facile
de raisonner sur la somme (\ref{equationroue}) qui représente le $\log$ formel de la fonction de
densité.

On calcule ensuite la différentielle en $\xi$ comme
précédemment.

\begin{figure}[!h]
\begin{center}
\includegraphics[width=5cm]{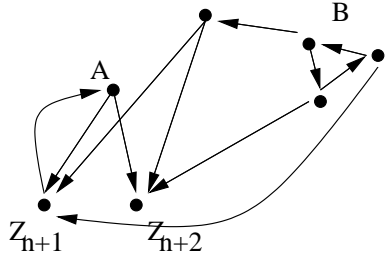}
\caption{\footnotesize Graphe Type de type roue}\label{ABroue}
\end{center}
\end{figure}

Seules importent les contractions d'arêtes qui arrivent sur les points représentant
$\xi$, notre paramètre de déformation. En effet les autres contributions se compensent par Jacobi
ou pour des raisons de dimension.
Les graphes que l'on traite sont simples et n'ont donc qu'une seule roue. Par
conséquent les graphes  typiques que l'on obtient par contraction des arêtes sont de deux types.
Le premier est comme dans la figure [\ref{ABroue}], avec $A$ un graphe simple de type Lie, une arête
joignant $z_{n+1}$ (resp. $z_{n+2}$) au sommet de $A$ et $B$ un graphe simple de type roue.
On note $w_{A,B}^{x}$ (resp. $w_{A,B}^{y}$) la $1$-forme obtenue par intégration sur les
configurations contenues dans $A$ et $B$. Comme dans la proposition \ref{produit} la partie
$1$-forme en $\xi$ ne peut provenir que des arêtes de $A$. On a donc encore une fois une situation
produit.

La deuxième possibilité consiste en un graphe comme le précédent mais sans le sous-graphe $B$.
C'est la cas précédent avec $B=\emptyset$.
Dans ce cas on note  $w_{A}^{x}$ (resp.$w_{A}^{y}$) la $1$-forme obtenue par intégration sur les
configurations contenues dans $A$.

Par conséquent les graphes simples qui dégénèrent en ces graphes typiques, sont comme précédemment
obtenus par adjonction d'une arête au point $z_{n+1}$ (resp. $z_{n+2}$) soit sur le graphe$B$
soit sur le graphe $A$. Toutefois si l'adjonction provient de $A$ et si $B$ n'est pas vide, le
graphe en question n'est pas simple; Par conséquent quand $B$ est non vide, l'adjonction provient de $B$. Lorsque $B$ est vide l'adjonction vient de $A$.
Par regroupement on obtient les contributions portant sur $x$ suivantes

$$\sum_{A, B}w_{A, B}^{x}(\xi)[x, A(x,y)]\cdot \partial_{x}B(x,y)=
\sum_{A, B}w_{A}^{x}(\xi)[x, A(x,y)]\cdot \partial_{x}w_{B}(\xi) B(x,y)=$$
\begin{equation}
[x, \Omega_{F}(x,y)]\cdot  \partial_{x}(\sum_{B}w_{B}(\xi)B(x,y))
 \end{equation}
et

\begin{equation}\label{equationroue2}
\sum_{A}w_{A}^{x}(\xi) \tr(\ad x \partial_{x}A (x,y)).
 \end{equation}

On a des contributions analogues portant sur $y$.

En conséquence on retrouve les $1$-formes du théorème \ref{deformationCBH}. Il vient maintenant
le théorème suivant.
\begin{theo}\label{deformationDensite}
La fonction de densité  $D(x,y)$ se déforme naturellement en la fonction $D_{\xi}(x,y)$ vérifiant
$$\mathrm{d}D_{\xi}(x,y)= \mathrm{tr}_{\g}(\mathrm{ad}x\partial_{x}\Omega_{F}(x,y)+
\mathrm{ad}y\partial_{y}\Omega_{G}(x,y)) D_{\xi}(x,y)+$$
$$[x, \Omega_{F}(x,y)]\cdot\partial_{x}D_{\xi}(x,y) +[y, \Omega_{G}(x,y)]
\cdot\partial_{y}D_{\xi}(x,y).$$
\end{theo}
{\it Preuve}:
On additionne les contributions précédentes. On remarque alors que l'on a raisonné sur le
$\log$ formel de $D_{\xi}(x,y)$, par conséquent la somme des contributions s'écrit

$$\mathrm{d}\log D_{\xi}(x,y)=\mathrm{tr}_{\g}(\mathrm{ad}x\partial_{x}\Omega_{F}(x,y)+
\mathrm{ad}y\partial_{y}\Omega_{G}(x,y))$$
$$[x, \Omega_{F}(x,y)]\cdot\partial_{x}\log D_{\xi}(x,y) +[y, \Omega_{G}
(x,y)]
\cdot\partial_{y}\log D_{\xi}(x,y).$$
C'est la formule cherchée car on a

$$ D_{\xi}(x,y) \mathrm{d}\log D_{\xi}(x,y)=\mathrm{d} D_{\xi}(x,y).$$
Comme on l'a dit dans l'introduction, nous travaillons au niveau des séries formelles mais
on peut justifier la convergence en $(x,y)$ dans un voisinage de l'origine de $0$ sans difficulté comme dans \cite{ADS}.$\Box$

\subsection{Conclusion et preuve du théorème principal}

Le théorème annoncé en introduction est maintenant démontré grace au théorème  \ref{deformationCBH} et au théorème \ref{deformationDensite}. Notre théorème montre alors, comme dans
\cite{KV}, que l'on dispose d'une déformation du produit de convolution qui assure que l'équation
(\ref{eqKV}) est vérifiée.

%
%

\end{document}